\documentclass[review]{elsarticle}

\usepackage{lineno,hyperref,amsmath,amssymb,amsfonts}
\modulolinenumbers[5]

\journal{Reliability Engineering \& System Safety}

\begin{document}

\begin{frontmatter}

\title{Uncertainty and sensitivity analysis of functional risk curves based on Gaussian processes}

\author{Bertrand Iooss and Lo\"ic Le Gratiet}
\address{EDF Lab Chatou, 6 Quai Watier, 78401 Chatou, France}

\begin{abstract}
A functional risk curve gives the probability of an undesirable event as a function of the value of a critical parameter of a considered physical system.
In several applicative situations, this curve is built using phenomenological numerical models which simulate complex physical phenomena.
To avoid cpu-time expensive numerical models, we propose to use Gaussian process regression to build functional risk curves.
An algorithm is given to provide confidence bounds due to this approximation.
Two methods of global sensitivity analysis of the models' random input parameters on the functional risk curve are also studied.
In particular, the PLI sensitivity indices allow to understand the effect of misjudgment on the input parameters' probability density functions.
\end{abstract}


\end{frontmatter}


\section{Introduction}

In industrial practice, the estimation of a functional risk curve (FRC) is often required as a quantitative measure of the safety of a system \cite{der12}.
A FRC (also called ``functional risk criterion'' or ``fragility curve'' in some application domains \cite{pas14}) gives the probability of an undesirable event as a function of the value of a critical parameter of a considered physical system.
The estimation of this curve sometimes relies on deterministic phenomenological computer models which simulate complex physical phenomena.
Uncertain input parameters of this computer code are modeled as random variables while one of its scalar output variables becomes the studied random variable of interest. 

As it is based on a probabilistic modeling of uncertain physical variables and their propagation through a numerical model, this problem can be directly related to the uncertainty management methodology in numerical simulation \cite{der12,baudut17}.
This methodology proposes a generic framework of modeling, calibrating, propagating and prioritizing uncertainty sources through a numerical model (or computer code).
Indeed, investigation of complex computer code experiments has remained an important challenge in all domains of science and technology, in order to carry on simulations as well as predictions, uncertainty analysis or sensitivity studies \cite{fanli06,smi14}.
In this paper, the physical numerical model $G$ is expressed as
\begin{equation}\label{eq:model}
Y = G(a,X) = G(a,X_1,\ldots,X_d)\;,
\end{equation}
with $a\in \mathbb{R}$ the critical input parameter of the FRC and $X \in {\mathbb{R}}^d$ the random input vector of dimension $d$ and $Y \in \mathbb{R}$ a scalar model output.
In this paper, parameters $a$ and $X$ are considered as independent. However the presented results hold when $a$ and $X$ are dependent and the distribution of $X$ conditionally to $a$ is known.

However, standard uncertainty treatment techniques require many model evaluations and a major algorithmic difficulty arises when the computer code under study is too time expensive to be directly used.
For cpu-time expensive models, one solution consists in replacing the numerical model by a mathematical approximation, called a response surface or a metamodel. 
Several statistical tools based on numerical design of experiments, efficient algorithms for uncertainty propagation and metamodeling concepts will then be useful \cite{fanli06}.
In this paper, Gaussian process regression \cite{raswil06} is used as a metamodeling technique (surrogate model of the computer code) and applied in the particular context of a FRC as a quantity of interest.
Numerous studies have shown that this model provides a powerful statistical framework to compute an efficient predictor of a deterministic computer code response \cite{marioo08,legmar17}.

Associated to the estimation of the quantity of interest after the uncertainty propagation step, the sensitivity analysis step is performed to determine those parameters that mostly influence the model response \cite{salcha00}.
In particular, global sensitivity analysis methods take into account the overall uncertainty ranges of the input parameter (see \cite{ioolem15} for a recent review).
Several works have focused on estimating global sensitivity indices (especially the variance-based ones, known as the Sobol' indices) using Gaussian process models \cite{oakoha04,legcan14,legmar17}.
In this paper, we present global sensitivity indices attached to the whole FRC (firstly introduced in \cite{legioo17} in a particular context), while showing how to develop them with a Gaussian process model.
We also focus on a recently developed method, the Perturbed-Law based sensitivity Indices (PLI) of \cite{lemser13}, which seems promising in terms of computational efficiency and interpretation.

Two examples of FRC largely used in some industrial safety practices are given in the next section.
The third section describes the Gaussian process way to model and estimate a FRC.
The fourth section develops the sensitivity indices adapted to FRC.

\section{Motivating examples and classical methods}

\subsection{Probability of detection curves}\label{sec:POD}

In several industries, the Probability Of Detection (POD) curve is a standard tool to evaluate the performance of Non Destructive Testing (NDT) procedures. 
The goal is to assess the inspection capability for the detection of flaws in the inspected structure. 
For instance, new aeronautic regulations require to perform appropriate damage tolerance assessments to address the potential for failure due to critical parts of material and manufacturing variability.
This imposes enhanced requirement on POD sizing.
Another example is the eddy current non destructive examination process which is used to ensure integrity of steam generator tubes in nuclear power plants \cite{maucos13,broleg15,legioo17}.
Experimental campaigns and simulated experiments are used in order to provide POD curves, which can be used to demonstrate the performance of inspection process to regulatory authorities.

In practice, the high costs of the implementation of experimental POD campaigns combined with the continuous increase in the complexity of configuration make them sometimes unaffordable. 
To overcome this problem, it is possible to resort to numerical simulation of NDT processes.
This approach is commonly called MAPOD for ``Model Assisted Probability of Detection'' \cite{tho08,cal12,legioo17}. 

More precisely, in the POD context, the problem is formulated as follows.
Given a threshold $s>0$, a flaw is considered to be detected when $Y > s$ with $Y$ the signal amplitude and $s$ a detection threshold. 
Therefore the one dimensional POD curve is denoted by: 
\begin{equation}\label{eq:PoD}
 \forall a>0 \quad \mbox{POD}(a)= P\left( Y > s \mid a \right) \;,
\end{equation}
where $a$ is the POD parameter of interest (for example the size flaw) and $X$ are the other input parameters.
When simulated experiments are used to build the POD (see for example \cite{legioo17}), $Y$ is a scalar output of a numerical model (see Eq. \ref{eq:model})), where $a$ is determined by its bounds and $X$ is a random vector defined by its joint probability density function (pdf).

In Sections 3 and 4 of this paper, we will illustrate the mathematical methods for FRC building and analysis on a case study on the eddy current examination of the wear induced on steam generator tubes by anti-vibration bars \cite{maucos13}.
For details, the reader is referred to \cite{broleg15,legioo17}.
The study involves the finite element computer code C3D \cite{maucos13a}.
The NDT process under study is the examination of anti-vibration bars (AVB) subject to wear of steam generator tubes inside a nuclear power plant \cite{maucos13}.
The computer code used simulates by finite elements the eddy-current propagation phenomena.
The critical input parameter $a$ is the maximum of the two flaw sizes while five other input parameters (vector $X$) of the computer code are also considered (with their associated pdf): E (pipe thickness), ebav$_1$ and ebav$_2$ (lengths of the gap between the AVB and the first and second flaws), h$_{11}$ and h$_{12}$ (first and second flaw heights).
The output $Y$ of the computer code is the signal amplitude, and the threshold of detection $s$ is fixed at a specific value.

\subsection{Seismic fragility curves}\label{sec:seismic}

The ``fragility curve'' is a popular functional risk criterion, commonly used in many engineering fields \cite{pas14}.
It describes the probability that the actual damage to a structure exceeds a damage threshold, when the structure is assigned to a specified load intensity.
For instance, the seismic fragility curve, which concerns systems subject to earthquakes, is of particular interest in nuclear safety studies (see \cite{zen10}).

In the case of seismic risk assessment the load is usually expressed as a scalar characteristic of a seismic signal, typically the horizontal Peak Ground Acceleration (PGA), which is a common choice in civil engineering.
Several parameters and phenomena, distinct from the PGA, also influence the load: for each PGA value, the occurrence of the damage event is random. 
The fragility curve may then be interpreted as the cumulative distribution function of the ``structural
capability'', i.e. the maximum load the structure under investigation can bear without damage \cite{zenbor14,damkel14}.
Fragility curves are useful tools in structural analysis as they provide more complete information than the usual failure probability (established for a reference value of the load only). 

In standard practice, the assessment is made either following an approach entirely based on available expertise or by the statistical analysis of actually observed or simulated data. 
As actual damage data may be scarce due to the rarity of severe earthquakes liable to generate damages on highly safe structures, observations are generated by mock-up or (most often) numerical experiments. 

In this context, a fragility function is expressed as
\begin{equation}\label{eq:SFC}
 \forall a >0 \quad \mbox{Frag}(a)= P\left( Y > s \mid a \right) \;,
\end{equation}
where $Y$ is the variable which characterizes the structural damage, $s=Y_{\mbox{\tiny max}}$ is the maximal capacity of the structure and $a$ is an intensity measure characterizing the ground motion severity ({\it e.g.} peak ground acceleration).
When the fragility curve is estimated by numerical simulation (see for example \cite{zen10}), $Y$ is the numerical model output (see Eq. (\ref{eq:model})), $a$ is determined by its bounds and the random vector $X$ is defined by its joint pdf.

\subsection{Classical methods}\label{sec:clas}

The standard practice for establishing a fragility curve consists in assuming a lognormal distribution for $Y$ \cite{shifen00}.
The classical approach for determining a POD curve is called the Berens model and relies on assuming a simple linear relation (with Gaussian residuals) between $\log(Y)$ and $a$ \cite{ber88}.
In fact, both cases are equivalent and lead to the following form for the FRC:
\begin{equation}\label{eq:PODclassic}
P\left( Y > s | a \right) = \Phi\left( \frac{a-\alpha}{\beta} \right) ,
\end{equation}
where $\Phi$ is the standard Gaussian distribution function, and $\alpha$ and $\beta$ are two constants to be estimated (by the maximum likelihood estimation for example) as functions of the $Y$-data sample and $s$.
The parameter $a$ can also be expressed in terms of a logarithm in this expression.
We note that if the Berens method is applied in a FRC framework (e.g. with $Y$ the maximal displacement of a building and $a$ the pick ground acceleration of a seismic signal) the resulting PoD curve is the classical FRC.

As shown in \cite{legioo17}, this model can easily be improved by considering a Box-Cox transformation \cite{boxcox64} instead of a simple logarithmic one on $Y$: $\tilde{Y} = \frac{Y^{\lambda}-1}{\lambda}$. 
$\lambda$ is determined by maximum likelihood as the real number that offers the best linear regression of $\tilde{Y}$ regarding the parameter $a$.
The model is now based on $\tilde{Y}$ and is defined as
\begin{equation}
\tilde{Y}(a) = \beta_0 + \beta_1 a + \epsilon,
\end{equation}
with $\epsilon$ the model error such as $\epsilon \sim \mathcal{N}\left(0,\sigma_{\epsilon}^2 \right)$. 
The maximum likelihood method provides the estimators $\hat{\beta_0}, \hat{\beta_1} \; \text{and} \; \hat{\sigma_{\epsilon}}$. 
Hence the model implies the following result: $\forall a>0, \quad \tilde{Y}(a) \sim \mathcal{N}\left(\hat{\beta_0} +\hat{\beta_1}a, \hat{\sigma_{\epsilon}}^2 \right)$. 
Then the value of the FRC can be estimated from $\hat{\beta_0}, \hat{\beta_1} \; \text{and} \; \hat{\sigma_{\epsilon}}$, as its confidence interval from the property of the maximum likelihood estimators.
Note that it is possible to relax the Gaussian hypothesis on $\epsilon$ by using non-parametric approaches \cite{maikon17,legioo17}

\section{Gaussian-process based functional risk curve estimation}

Classical methods described in Section \ref{sec:clas} use a simple model between $Y$ and $a$ without modeling the functional dependence of $Y$ on the other uncertain variables $X$.
As in our case data $y_N$ (sample of $Y$ of size $N$) are obtained from numerical models, the input variables $X$ are controlled and can be introduced in the FRC determination process.
However, with cpu-time expensive numerical models, model evaluations can be somewhat limited and only small samples of $Y$ can be obtained.
As shown in \cite{demjen12,zenbor14,legioo17}, the use of a metamodel, which is a mathematical approximation of the computer code \cite{fanli06}, is useful.
We highlight that it is important in a safety framework to evaluate accurately the uncertainty on the surrogate model.
We propose in this section a general formulation for FRC that uses Gaussian process regression as a metamodel \cite{sacwel89,raswil06}.
We propose in this section a general formulation for FRC that uses the Gaussian process regression metamodel \cite{sacwel89,raswil06}. 
Through the Gaussian assumption, this metamodel provides a convenient framework to perform uncertainty quantification on the quantity of interest estimate.
This formulation has been introduced in \cite{legioo17} without the full mathematical algorithm which is developed in the present paper.

Let $X$ be the random vector of influential and uncertain parameters of the computer model $G(\cdot)$ and $a$ be the parameter of interest (the abscissa of the FRC).
The prior knowledge on $G(a,X)$ is modeled by $Y(a,X)$ and defined as follows .
\begin{equation} 
Y(a,X) = \beta_0 + \beta_1 a + Z_{\sigma^2,\theta}(a,X) , 
\end{equation}
where $Z_{\sigma^2,\theta}$ is a centered Gaussian process. We make the assumption that $Z_{\sigma^2,\theta}$ is second order  stationary with variance $\sigma^2$ and a parametric covariance function (which corresponds to the kernel in other statistical learning methods).
As covariance choice, one can cite the popular Mat\'ern 5/2 parametrized by its lengthscale $\theta$).
Thanks to the maximum likelihood method, we can estimate the values of the so far-unknown parameters: $\beta_0, \beta_1, \sigma^2$ and $\theta$ (see for instance \cite{marioo08} for more details).

Gaussian process regression (also known as the kriging process) provides an estimator of $G(a,x)$ which is called the kriging predictor and written $\widehat{Y}(a,x)$.
In addition to the kriging predictor, the kriging variance $\sigma_Y^2(a,x)$ quantifies the uncertainty induced by estimating $Y(a,x)$ with $\widehat{Y}(a,x)$. 
The predictive distribution is given by $Y(a,x)$ conditioned by $y_N$:
\begin{equation}\label{eq:kriging}
\forall x \quad \left(Y(a,x) \mid y_N \right) = {\mathcal{N}} \left( \widehat{Y}(a,x), \sigma_Y^2(a,x) \right) 
\end{equation}
where $\widehat{Y}(a,x)$ (the kriging mean) and $\sigma_Y^2(a,x)$ (the kriging variance) can both be explicitly estimated. 

Obtaining the FRC $\Psi(a)$ consists in replacing $Y = G(a,X)$ by its Gaussian process metamodel (\ref{eq:kriging}) in (\ref{eq:PoD}).
 Hence we can estimate the value of $\Psi(a)$, for $a>0$ from:
\begin{equation}\label{eq:PODkriging}
\Psi(a) = P \left(Y_N(a,X) > s | a \right) ,
\end{equation}
where $Y_N(a,X)$ is a Gaussian process which follows the distribution $(Y(a,X) \mid y_N)$ and $s$ is a fixed threshold value.
Two sources of uncertainty have to be taken into account in (\ref{eq:PODkriging}): the first coming from the parameter $X$ and the second coming from the Gaussian distribution in (\ref{eq:kriging}).

From (\ref{eq:PODkriging}), the following estimate for $\Psi(a)$ can be deduced:
\begin{equation}\label{eq:Psi} 
\widehat \Psi(a) = E_X \left[ 1 - \Phi\left(\frac{s-\widehat{Y}(a,X)}{\sigma_Y(a,X)}\right)\right].
\end{equation}
This equation corresponds to the mean FRC with respect to the Gaussian process metamodel.
The expectation in (\ref{eq:Psi}) is estimated using a classical Monte Carlo integration procedure.

In order to estimate the uncertainty on the FRC estimation, we start from its integral expression:
\begin{equation}
\Psi(a) = P_X (G(a,X) > s |a) = \int{1_{G(a,x) > s}} \, f(x) d{x},
\end{equation}
where $f(x)$ is the joint pdf of $X$ (independent on $a$).
The first uncertainty source on $\Psi(a)$ comes from the numerical evaluation of the integral inside the $\Psi(a)$ definition.
This evaluation is done by the Monte Carlo method:
\begin{equation}
\Psi(a) \approx \Psi_{MC}(a) =  \frac{1}{n} \sum_{i=1}^{n} 1_{G(a,x^{(i)}) > s},
\end{equation}
where $(x^{(i)})_{i=1,\dots,n}$ is a sample of the random variable $X$. 
The required size of this sample can be too large to use this formula in practice (case of a costly evaluation of the code).
We then replace the code by its Gaussian process approximation $Y_N(a,X)$:
\begin{equation}
\Psi_{MC}(a) \approx \Psi_{MC,GP}(a) = \frac{1}{n} \sum_{i=1}^{n} 1_{{Y_N(a,X)}  > s}.
\end{equation}
Let us recall that $Y_N(a,X)$ is a Gaussian process with known mean and variance.
To compute the integration error, the central limit theorem is used and gives for $n \rightarrow \infty$:
\begin{equation}\label{eq:tcl}
\sqrt{n} \left(\Psi_{MC,GP}(a) - \Psi_{GP}(a)\right) \longrightarrow {\mathcal{N}}\left(0,\Psi_{GP}(a)(1-\Psi_{GP}(a))  \right),
\end{equation}
where
\begin{equation}
\Psi_{GP}(a) = \int{1_{{Y_N(a,X)} > s}} \, f(x) d{x}.
\end{equation}
$\Psi_{GP}(a)$ corresponds to the FRC for the process $Y_N(a,X)$.
The algorithm uses a significantly large value of $n$ in order to deal with a valid Gaussian approximation.

The second uncertainty source on $\Psi(a)$ estimation comes from Gaussian process approximation.
We simulate $m$ realizations $(y^{(j)}(a,x))_{j=1,\dots,m}$ from $({Y(a,x) | y_N})$ in order to evaluate the variability of $\Psi_{MC,GP}(a)$ which comes from Gaussian process approximation.
Thus, we compute:
\begin{equation}
 \psi_{MC,GP}^{(j)}(a) = \frac{1}{n} \sum_{i=1}^{n} 1_{{y^{(j)}(a,x^{(i)})} > s}.
\end{equation}
Then, for each $\psi_{MC,GP}^{(j)}(a)$, $j=1,\dots,m$, we compute the Monte Carlo error using the central limit theorem cited before in (\ref{eq:tcl}) where the variance of the Gaussian distribution in \eqref{eq:tcl} is estimated with $\psi_{MC,GP}^{(j)}(a) (1-\psi_{MC,GP}^{(j)}(a) )$, $j=1,\dots,m$.

From \eqref{eq:tcl}, we compute for each $j=1,\dots,m$ a sample of size $n_{CLT}$.
From this double Monte Carlo method,  we obtain a sample $\left[\psi_{MC,GP}^{(j,k)}(a)\right]_{\stackrel{j=1,\dots,m }{k=1,\dots,n_{CLT}}}$ of size $m \times n_{CLT}$. Then the FRC estimate is given by :
\begin{displaymath}
\frac{1}{m \times n_{CLT}} \sum_{\stackrel{j=1,\dots,m }{k=1,\dots,n_{CLT}}} \psi_{MC,GP}^{(j,k)}(a),
\end{displaymath}
and the uncertainty is evaluate from  $\left[\psi_{MC,GP}^{(j,k)}(a)\right]_{\stackrel{j=1,\dots,m }{k=1,\dots,n_{CLT}}}$ by using classical quantile estimators.

\begin{figure}[!ht]
  \centering
	\includegraphics[scale=0.6]{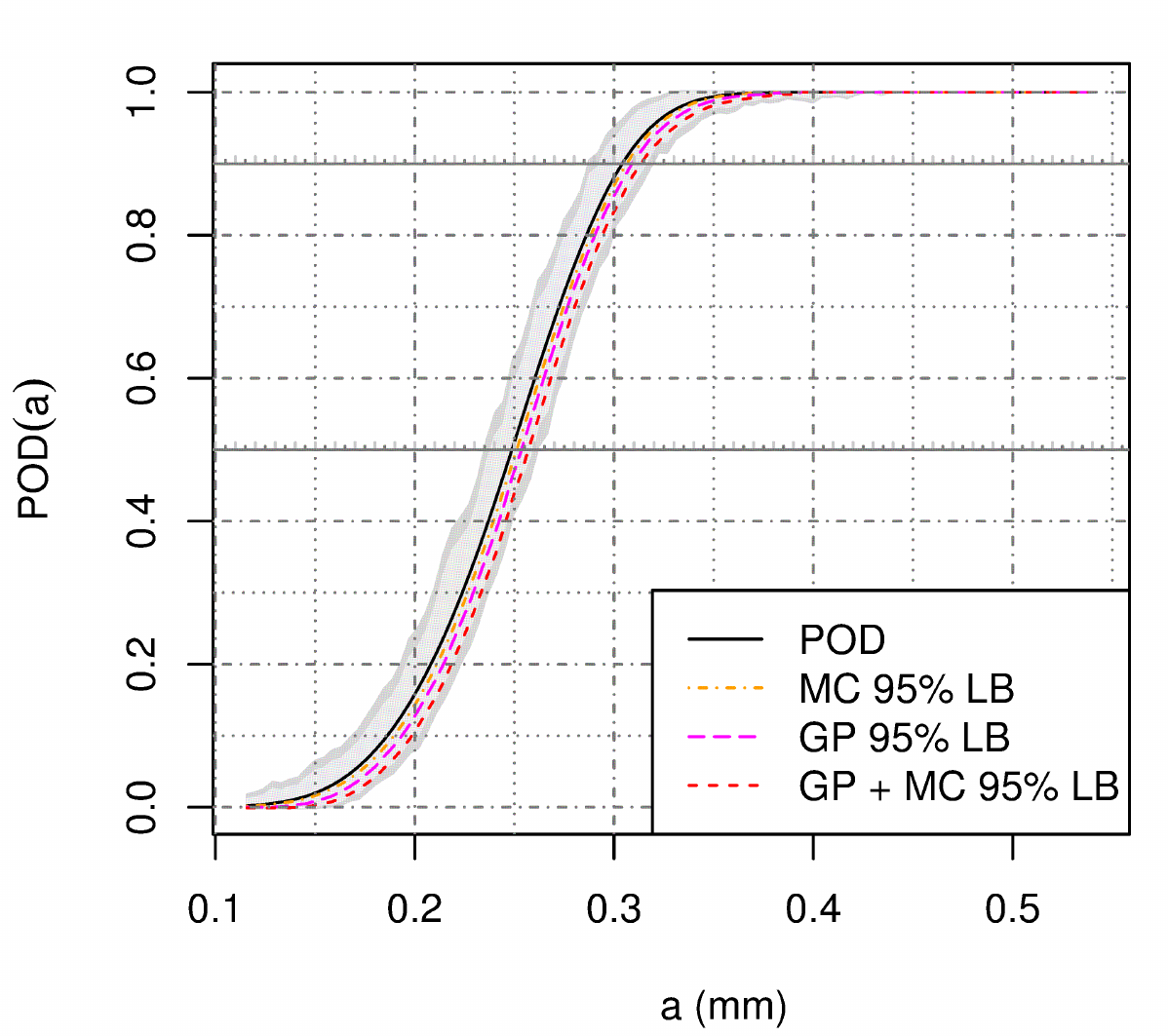}
\caption{Example of FRC estimated with a Gaussian process model (from \cite{legioo17}). The present curve is a POD curve. The lines represent unilateral confidence intervals (LB means ``lower bound''). The shaded areas represent the minimal and the maximal observations of POD($a$) when only the uncertainty on the Gaussian process is considered (light gray) and when both the Gaussian process and the Monte-Carlo errors are considered (dark gray).}
  \label{fig:krigingPoD}
\end{figure}

From this double Monte Carlo method, we are able to compute the FRC uncertainty due to the Gaussian process metamodel and due to the numerical integration process.
Figure \ref{fig:krigingPoD} illustrates the algorithm on a POD curve estimation in an NDT application (see section \ref{sec:POD}).
A Gaussian process metamodel has been built on $N=100$ numerical simulations of eddy current examinations of steam generator tubes \cite{legioo17}.
We visualize the confidence interval induced by the Monte Carlo (MC) estimation ($n = 10,000$, $n_{CLT} = 100,000$), the one induced by the Gaussian process approximation ($m=3,000$) and the total confidence interval (including both approximations: GP+MC). It could be noticed that from the double Monte-Carlo sample $\left[\psi_{MC,GP}^{(j,k)}(a)\right]_{\stackrel{j=1,\dots,m }{k=1,\dots,n_{CLT}}}$, one can evaluate with an ANOVA decomposition the contribution of the Gaussian process and the parameters $X$ on the FRC estimate mean squared error (see \cite{legcan14}).

Choosing $m$, the number of conditional Gaussian process realizations to simulate, is a difficult issue of this process. 
As explained in \cite{legcan14} and other papers, there are some computational constraints (in memory and cpu time) when simulating conditional Gaussian processes on large samples. 
The choice of $m=3000$ has been chosen as a sufficient value (that we avoid to be too large) in order to estimate a $95\%$ lower bound (due to the Gaussian process error). 
Then $n$ has been tuned in order to have a smaller error induced by the Monte Carlo sampling process than the one induces by the Gaussian process. 

\section{Global sensitivity analysis of functional risk curves}

The objective of sensitivity analysis is to determine those input variables that mostly influence the model response \cite{salcha00,borpli16,ioosal17}.
Global sensitivity analysis methods take into account the overall uncertainty of the input parameters.
Previous works on the POD \cite{rupbla14} and seismic fragility context \cite{borzen13} have only considered sensitivity analysis on the variance and distribution of the model output variable $Y$.
In this section, we propose two global sensitivity indices attached to the whole FRC: one already introduced in \cite{legioo17} (aggregated Sobol' indices) and a new one based on a completely different idea (perturbed-law based indices).
Therefore, the influence of the parameter of interest $a$, as it is the FRC abscissa, will not be considered as an input variable whose sensitivity index has to be computed.

\subsection{Aggregated Sobol' indices}

We first focus on the variance-based sensitivity formulation \cite{sob93}, which is one the most popular tools that has been proven to be robust and easily interpretable \cite{ioolem15}.
The associated sensitivity indices are often called Sobol' indices.

In the case of independent inputs, using the Sobol-Hoeffding decomposition \cite{pritar17}, the variance of the numerical model $Y=G(X_1,\ldots,X_d)$ is decomposed into the following sum:
\begin{equation}\label{eq:Var}
\mbox{Var}(Y) = V = \sum_{i=1}^d V_i + \sum_{1 \leq i < j \leq d} V_{ij} + \ldots + V_{1\ldots d} \,
\end{equation}
with $V_i=\mbox{Var}[E(Y|X_i)]$, $V_{ij}=\mbox{Var}[E(Y|X_i, X_j)]- V_i - V_j$, etc.
Then, $ \forall i,j=1\ldots d, \;i<j$, the Sobol' indices are given by
\begin{equation}\label{eq:Sobol}
S_i = \frac{V_i}{V} \;, S_{ij} = \frac{V_{ij}}{V} \;, \ldots, \mbox{ and } T_i = S_i + \sum_{j \neq i} S_{ij}+\sum_{j<k,j \neq i, k \neq i} S_{ijk} + \ldots \;.
\end{equation}
The first-order Sobol' index $S_i$ measures the individual effect of the input $X_i$ on the variance of the output $Y$, while the total Sobol' index $T_i$ measures the $X_i$ effect and all the interaction effects between $X_i$ and the other inputs (as the second-order effect $S_{ij}$).
$T_i$ can be reinterpreted as $T_i=\displaystyle 1 - \frac{V_{-i}}{V}$ with $V_{-i}=\mbox{Var}[E(Y|X_{-i})]$ and $X_{-i}$ is the vector of all inputs except $X_i$.
Each Sobol' index is then interpreted in terms of percentage of total variance explanation.

In order to define similar sensitivity indices for the whole FRC (Eq. (\ref{eq:Psi})), \cite{legioo17} defines the following quantities:
\begin{equation}\label{eq:PoDX}
\begin{array}{rcl}
\Psi_X(a) &= &P_X( Y > s \mid a,X) = 1 - \Phi\left(\frac{s-\widehat{Y}(a,X)}{\sigma_Y(a,X)}\right) \;,\\
\Psi_{X_i}(a) &= &P_X( Y > s \mid a,X_i) \;,\\
\Psi_{X_{-i}}(a) &= &P_X( Y > s \mid a,X_{-i}) \;,\\
D & = &E_X\|\Psi(a)-\Psi_X(a)\|^2\\
 & = & E_X \left[ \int \left(\Psi(a)-\Psi_X(a) \right)^2\,da \right],
\end{array}
\end{equation}
with $\|.\|$ the euclidean norm.
As in \cite{gamjan14} (see also \cite{marsai17}) which deals with functional model outputs, our objective is to obtain a single sensitivity index for each input by averaging the variability of the function, here the FRC depending on $a$.
The specificity here is that we deal with a FRC function.
Therefore, \cite{legioo17} defines the FRC aggregated Sobol' indices by:
\begin{equation}
\begin{array}{rcl}
S_i^{\mbox{\tiny FRC}} & = & \displaystyle \frac{E_X \|\Psi(a)-\Psi_{X_i}(a)\|^2}{D}\;,\\
T_i^{\mbox{\tiny FRC}} & = & \displaystyle \frac{E_X \|\Psi_X(a)-\Psi_{X_{-i}}(a)\|^2}{D}\;.
\end{array}
\end{equation}
$S_i^{\mbox{\tiny FRC}}$ (resp. $T_i^{\mbox{\tiny FRC}}$) gives the first-order (resp. total) effect of $X_i$ on the mean FRC.

These FRC aggregated Sobol' indices are easily computed by the metamodels.
In particular, the kriging metamodel allows to replace $P_X( Y > s \mid a)$ by the expectation $\displaystyle E_X \left[ 1 - \Phi\left(\frac{s-\widehat{Y}(a,X)}{\sigma_Y(a,X)}\right)\right]$ in the FRC expressions of (\ref{eq:PoDX}).
In order to compute the Sobol' indices, we use the classical pick-freeze formulas (see for example \cite{pritar17}).
The confidence intervals due to the finite Monte Carlo sample are then obtained by bootstrap.
Let us noticed that it is also possible to estimate the error on the aggregated Sobol' indices due to the kriging error as in \cite{legcan14} for the classical Sobol' indices.
This computationally heavy process has not been applied here.

\begin{figure}[!ht]
  \centering
	\includegraphics[width=0.49\textwidth]{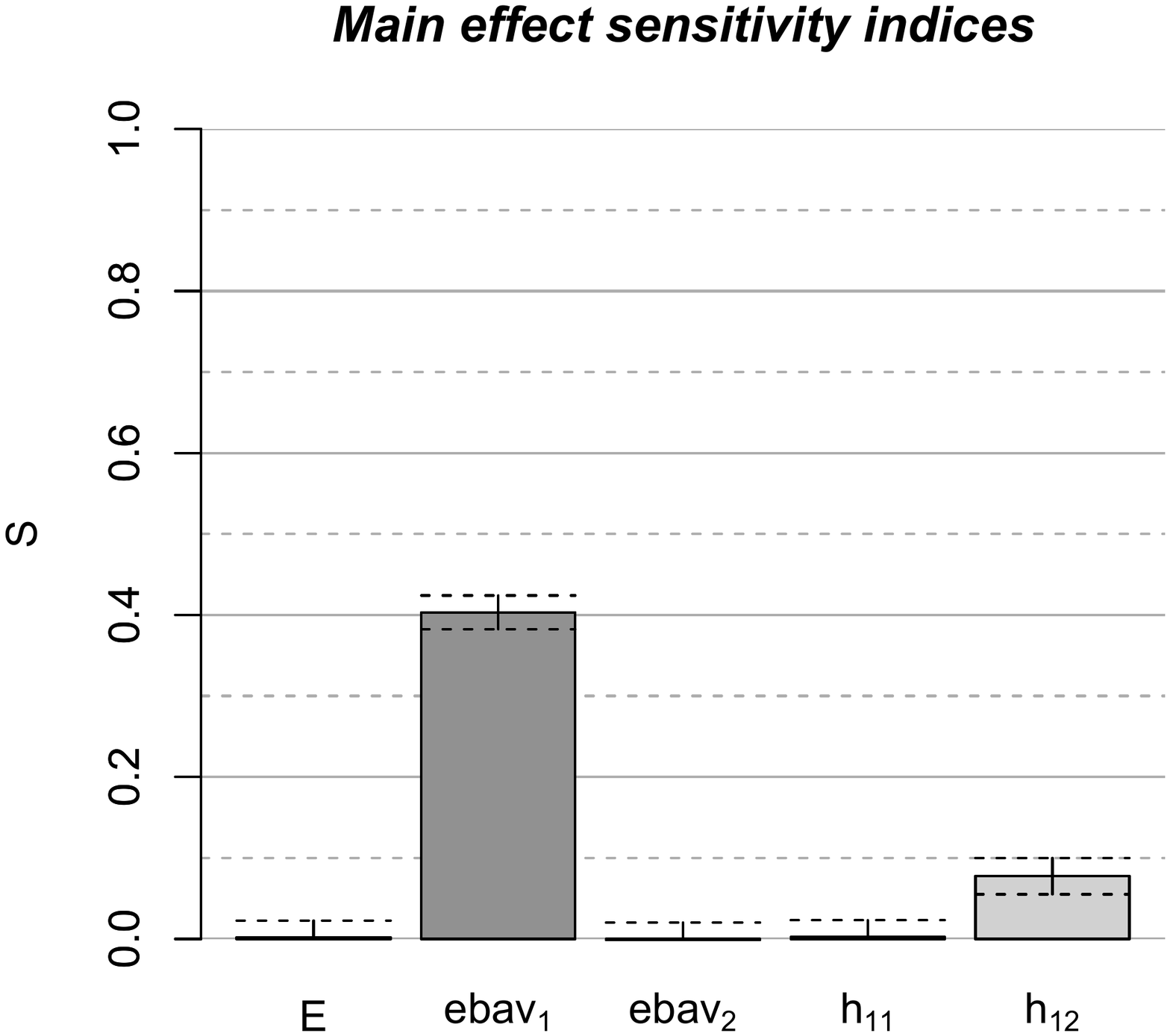}
	\includegraphics[width=0.49\textwidth]{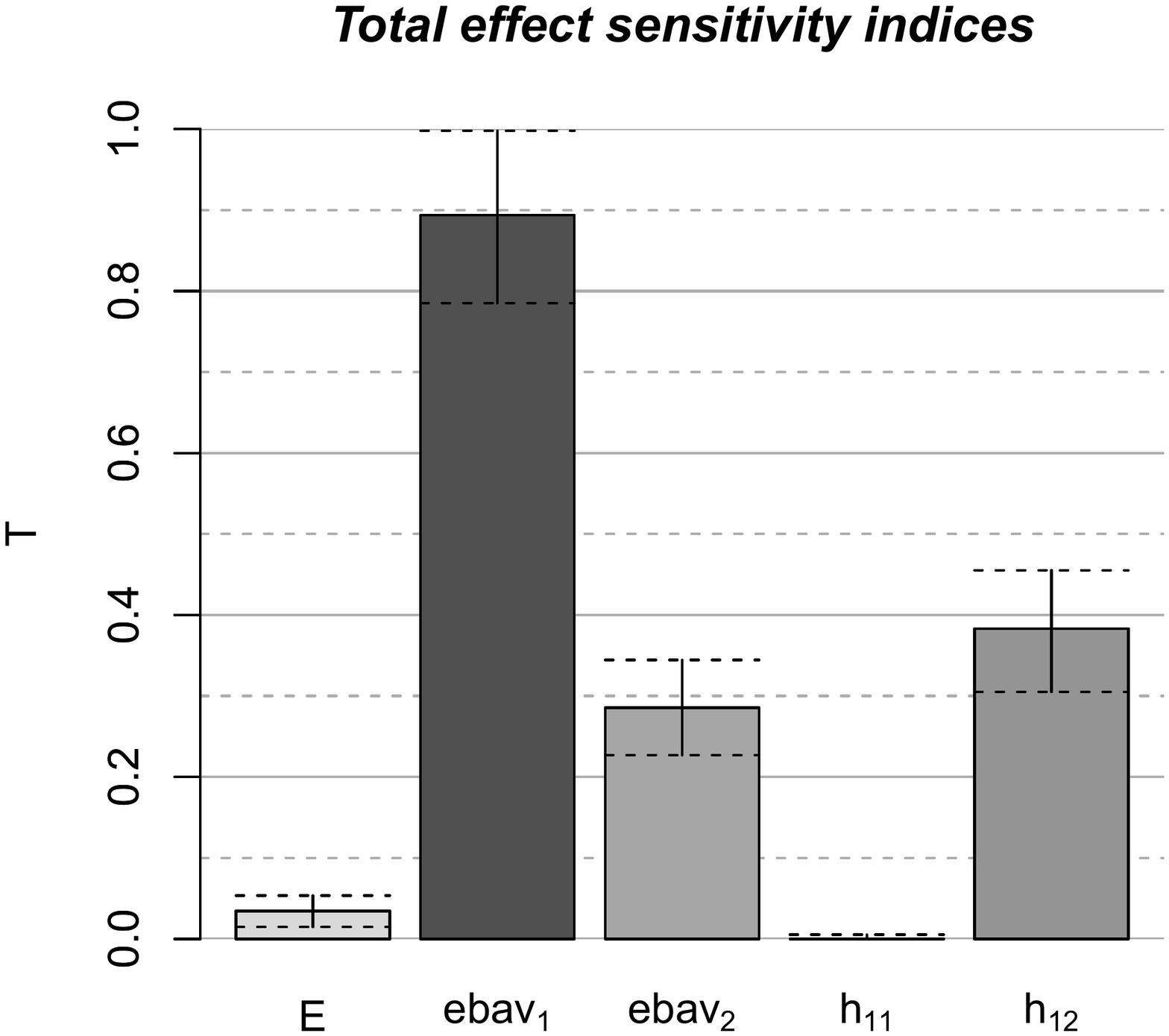}
	
	\vspace{-2cm}
\caption{First order (left) and  total (right) FRC aggregated Sobol' indices (from \cite{legioo17}).}
  \label{fig:SobolPoD1}
\end{figure}

The sensitivity analysis results on our use case are given in Figure \ref{fig:SobolPoD1} with the FRC aggregated Sobol' indices of the five physical input parameters (called E, ebav$_1$, ebav$_2$, h$_{11}$ and h$_{12}$).
ebav$_1$ is the main influent parameter on the FRC (the POD curve of Figure \ref{fig:krigingPoD}), followed by h$_{12}$ and ebav$_2$.
$E$ and h$_{11}$ have no influence.

\begin{figure}[!ht]
  \centering
	\includegraphics[width=0.49\textwidth]{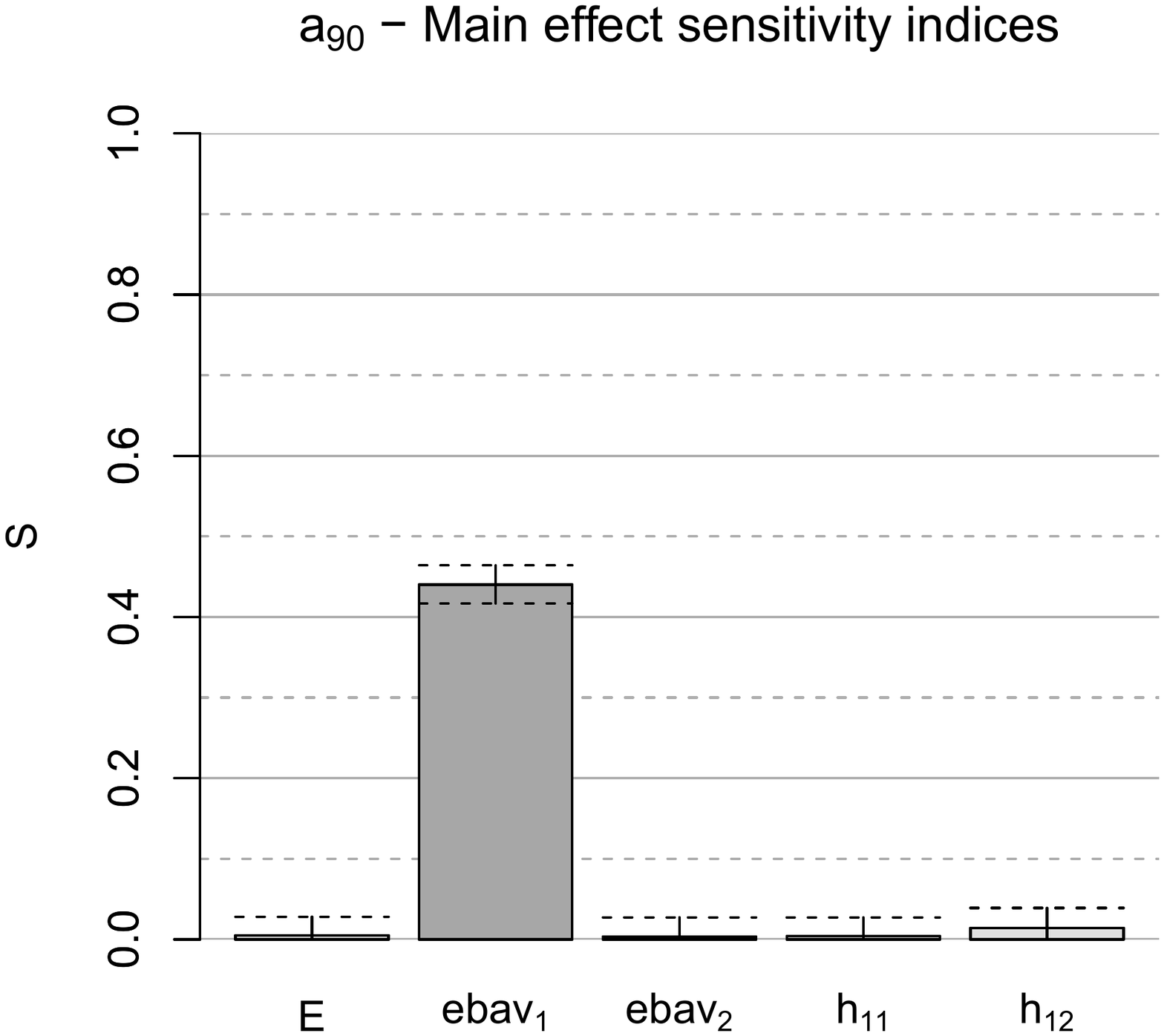}
	\includegraphics[width=0.49\textwidth]{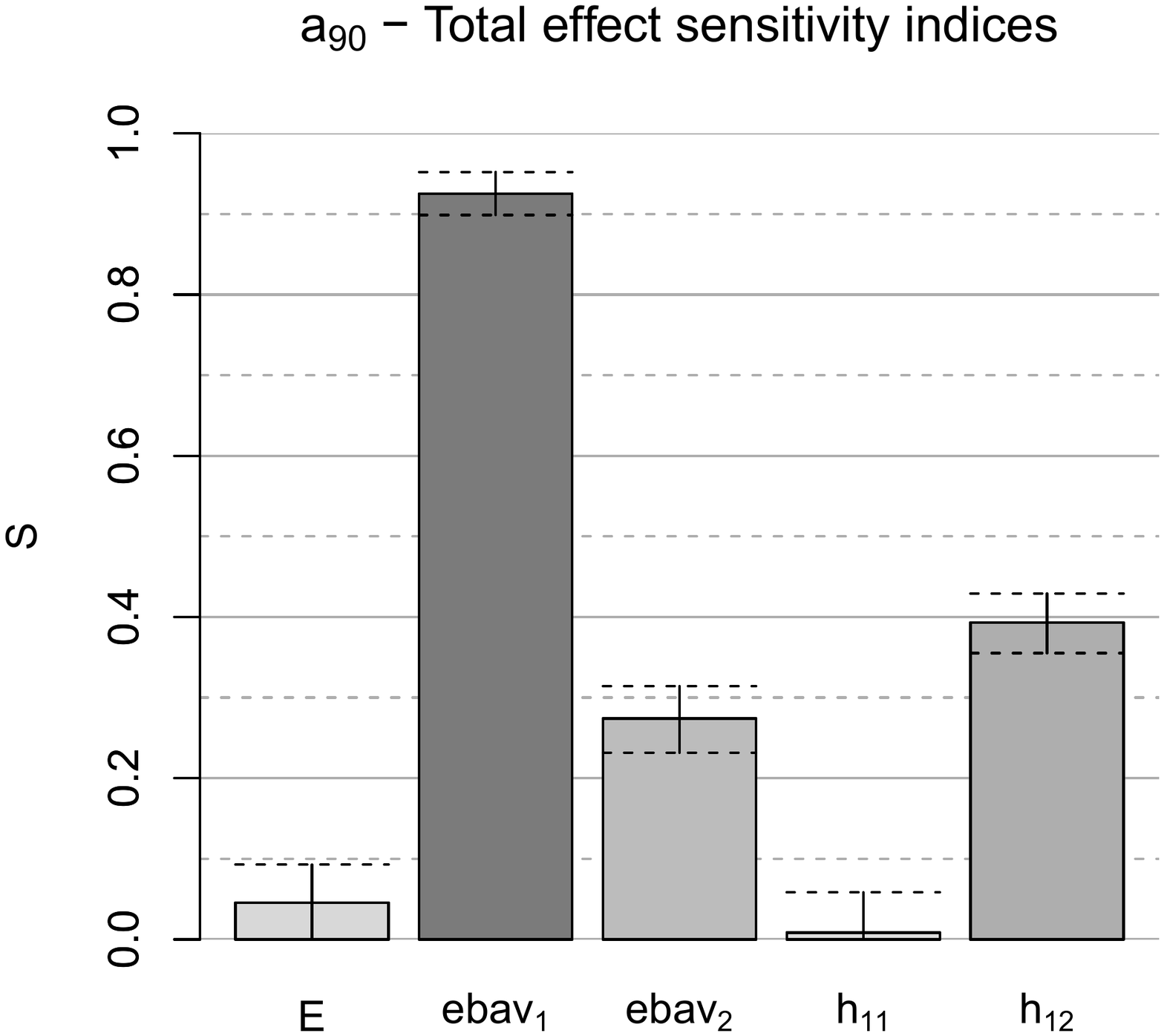}
	
	\vspace{-2cm}
\caption{First order (left) and total (right) Sobol' indices on $a_{90}$ (from \cite{legioo17}).}
  \label{fig:SobolPoD2}
\end{figure}

The FRC aggregated Sobol' indices quantify the sensitivity of each input on the overall curve.
However, we could be interested in the sensitivities on the FRC at a specific value of the parameter of interest $a$.
This can be directly done by replacing $Y$ by $\Psi_X(a)$ in the equations (\ref{eq:Var}) and (\ref{eq:Sobol}).
If we are now interested by the sensitivities at a specific probability value of the FRC, we have to study the inverse function of the FRC: $\Psi_X^{-1}(p)$ with $p$ a given probability.
Similarly to the previous case, the Sobol' indices can be obtained by replacing $Y$ by $\Psi_X^{-1}(p)$ in the equations (\ref{eq:Var}) and (\ref{eq:Sobol}).
Figure \ref{fig:SobolPoD2} gives these sensitivity indices on our data for $p=0.90$ which corresponds to the quantile of $a$ at order $90\%$ (noted $a_{90}$).
We see that the influences of the inputs on $a_{90}$ are very close to the ones on the overall POD curve.
It could be explained by the fact that the different inputs have rather linear effects on the output of the model (the amplitude).

\subsection{Perturbed-law based sensitivity indices}
 
We propose now to quantify the impact on the FRC of a perturbation of the input parameters pdfs by answering to the following question: ``What would be the FRC if the pdf of the $i^{\mbox{\tiny th}}$ input $X_i$ had been modified?''.
In this approach, all the input parameters are modeled by random variables, and their input probability densities are supposed to be unknown.
Let us remark that a negligible sensitivity index of an input will not allow to fix this input, but just to say that its pdf has no influence on the FRC.

We use the so-called Perturbed-Law based sensitivity Indices (PLI) measures recently introduced in \cite{lemser13} (see also \cite{borioo17}).
We start from the integral-form of the FRC:
\begin{displaymath}
P_X \left( G(a,x) > s \right) = \int 1_{G(a,x) > s} f(x)\,dx,
\end{displaymath}
where $f(x)$ is the joint pdf of $X$. 
Modifying the pdf $f_i(x_i)$ of $X_i$ gives us $f_{i,\delta}(x_i)$, the perturbed pdf of $X_i$.
After this perturbation, the FRC, denoted $P_{X_i,\delta}(\cdot)$ instead of $P_X(\cdot)$, can be written as:
\begin{equation}
P_{X_i,\delta} \left( G(a,x) > s \right) = \int { 1_{G(a,x) > s} \frac{f_{i,\delta}(x_i)}{f_i(x_i)} f(x)\,dx,}
\end{equation}

The PLI measures only consist of the comparison of the FRCs before and after the perturbation, and are defined by:
\begin{equation}
S_{i,\delta} = \left\{
\begin{array}{l}
\displaystyle \frac{ P_{X_i,\delta} \left( G(a,X) > s \right) - P_X \left( G(a,X) > s \right)  }{  P_{X_i,\delta} \left( G(a,X) > s \right)}  \\
\mbox{if} \quad P_{X_i,\delta} \left( G(a,X) > s \right) \geq  P_X \left( G(a,X) > s \right), \\
\\
\displaystyle \frac{ P_{X_i,\delta} \left( G(a,X) > s \right) - P_X \left( G(a,X) > s \right)  }{   P_X \left( G(a,X) > s \right) }\\
   \mbox{if} \quad P_{X_i,\delta} \left( G(a,X) > s \right) < P_X \left( G(a,X) > s \right). \\
\end{array}
\right. 
\end{equation}
A negative $S_{i,\delta}$ means that the FRC is smaller after the perturbation, while a positive $S_{i,\delta}$ means that the FRC has increased.
The estimation of $P_{X_i,\delta} \left( G(a,X) > s \right)$ is based on reverse importance sampling \cite{hes96}.
Asymptotical properties of the estimators give also some confidence intervals on the PLI measures. 

In \cite{lemser13}, the numerical model is directly used to estimate the PLI measures by large Monte Carlo samples.
We propose here to estimate the PLI measures using the Gaussian process metamodel and integrating its error in the estimates.
The mean FRC which we consider is given by (\ref{eq:Psi}):
\begin{equation}
\Psi(a) = \displaystyle E_X \left[ 1 - \Phi \left(\frac{s - \widehat{Y}(a,X) }{\sigma_Y(a,X)}\right) \right] = \displaystyle \int{ 1 - \Phi \left(\frac{s - \widehat{Y}(a,x) }{\sigma_Y(a,x)}\right)   f(x)\,dx}.
\end{equation}
After the perturbation of the pdf of the  $i^{\mbox{\tiny th}}$  input, the mean FRC is given by:
\begin{equation}
\Psi_{i,\delta}(a) = 
  \int{ 1 - \Phi \left(\frac{s - \widehat{Y}(a,x) }{\sigma_Y(a,x)}\right)  \frac{f_{i,\delta}(x_i)}{f_i(x_i)} f(x)\,dx}.
\end{equation}
PLI measures are then given by:
\begin{equation}
S_{i,\delta}(a) = \left\{
\begin{array}{l}
\displaystyle \frac{ \Psi_{i,\delta}(a)- \Psi(a)}{\Psi_{i,\delta}(a)} \quad \mbox{if} \quad \Psi_{i,\delta}(a) \geq \Psi(a), \\
\\
\displaystyle \frac{\Psi_{i,\delta}(a)-\Psi(a)  }{ \Psi(a) } \quad \mbox{if} \quad \Psi_{i,\delta}(a) <\Psi(a). \\
\end{array}
\right. 
\end{equation}

The last element of the PLI method is the definition of the perturbations  $f_{i,\delta}(x_i)$. 
\cite{lemser13} choose to perturb a statistical characteristic (for example the mean, or the variance, or a quantile, \ldots) of $X_i$ in order to be ``as close as possible'' to the initial pdf $f_{i}(x_i)$. 
The dissimilarity measure between $f_{i,\delta}(x_i)$ and $f_{i}(x_i)$, which contains the required properties, is the  Kullback-Leibler divergence:
 \begin{equation}
 KL(f_{i,\delta}(x_i),f_{i}(x_i)) = \int_{-\infty}^{+\infty}  f_{i,\delta}(x_i) \log \frac{f_{i,\delta}(x_i)}{f_{i}(x_i)}dx_i.
 \end{equation}
It implies that $f_{i,\delta}(X_i)$ and $f_i(X)$ have the same definition domain (as a consequence, the input domain bounds cannot be changed).
The Kullback-Leibler divergence has the advantage to be easily minimized \cite{lem14}, in order to obtain explicit solutions for $f_{i,\delta}(X_i)$ for classical pdf ({\it e.g.} Gaussian).

\begin{figure}[!ht]
  \centering
\vspace*{-3ex}	
	\includegraphics[trim=0 0 0 0,width=0.49\textwidth]{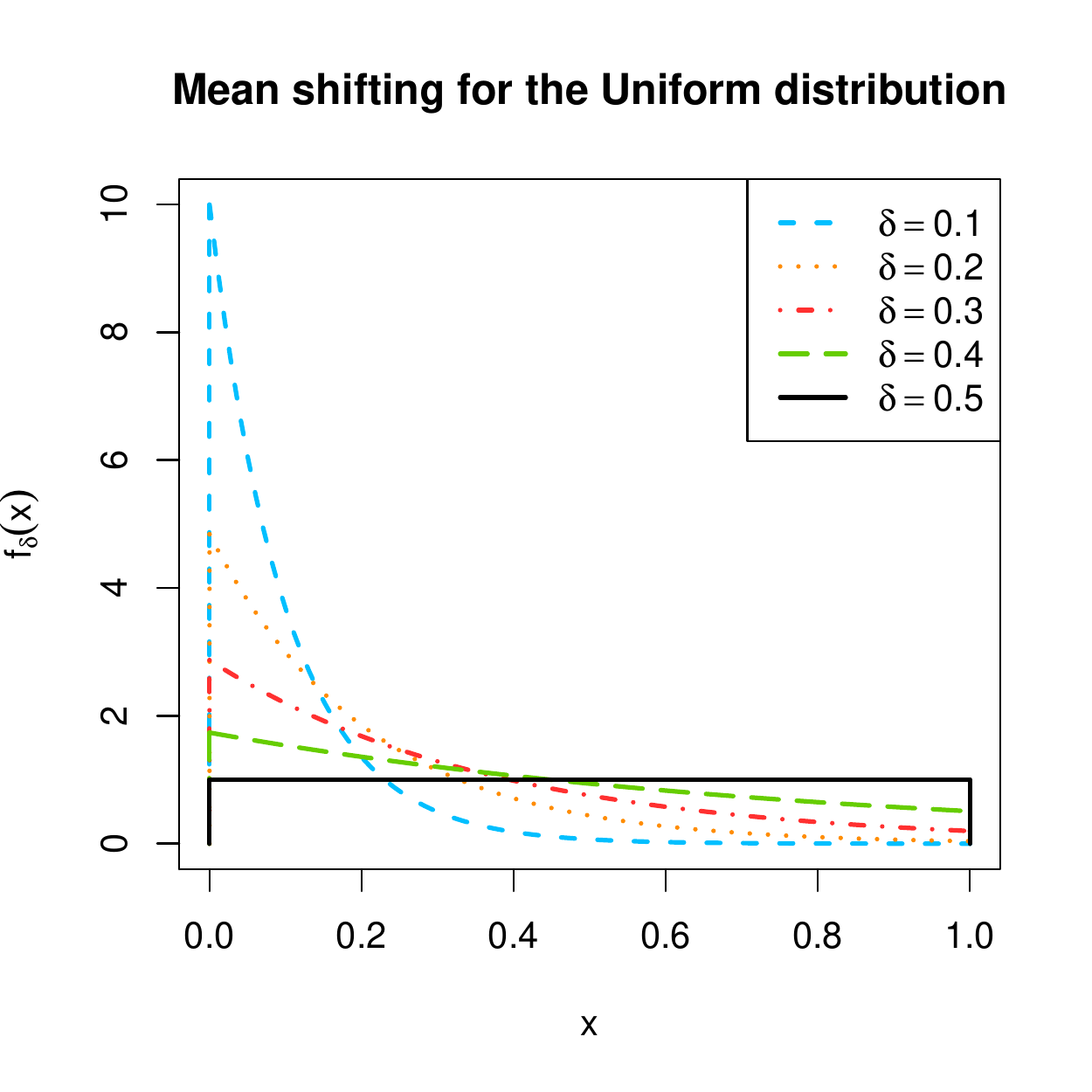}
		\includegraphics[trim=0 0 0 0,width=0.49\textwidth]{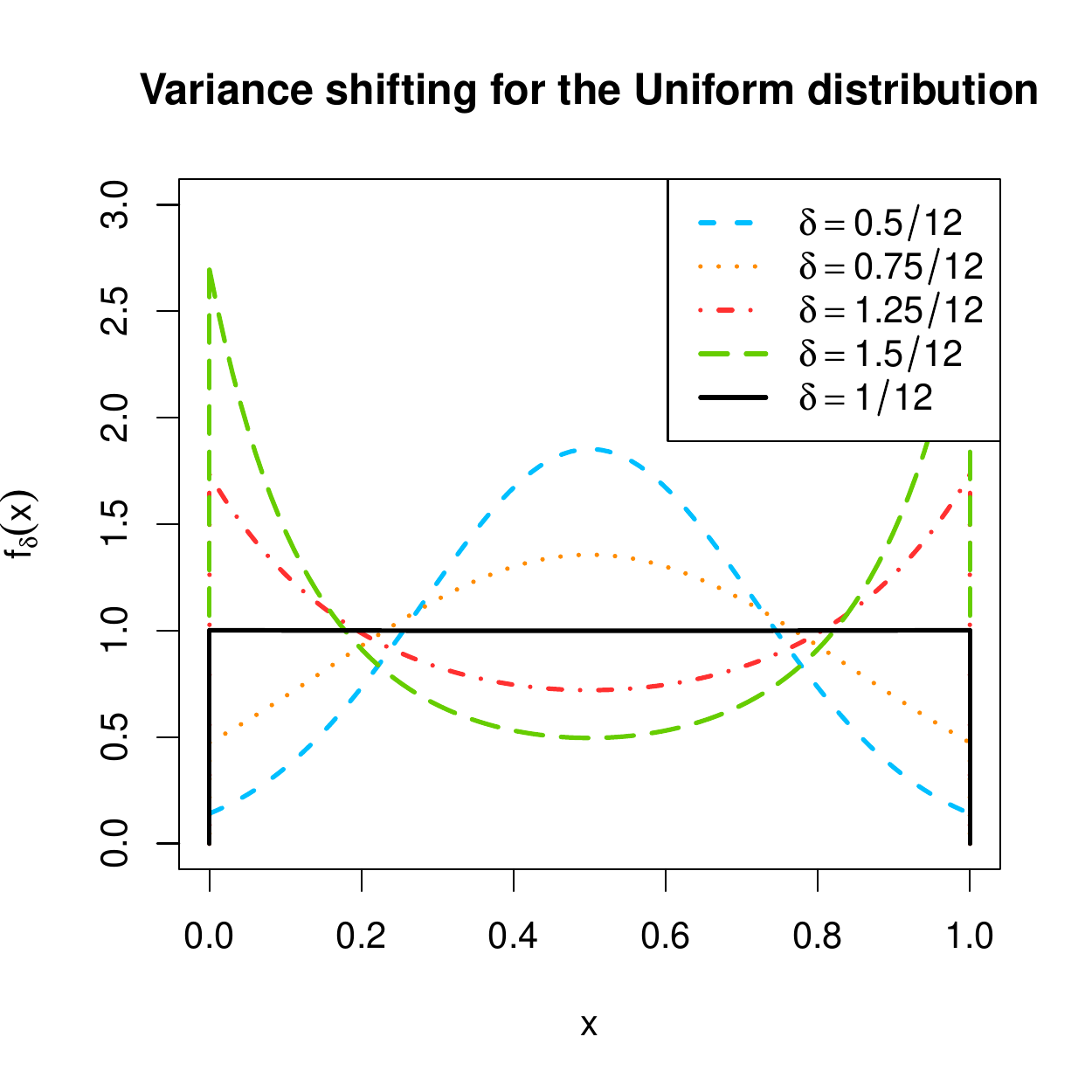}
  \caption{Examples of perturbations for the uniform pdf on $[0,1]$.
	Left: the mean of the pdf is perturbed ($\delta$ is the new mean value). 
	Right: the variance of the pdf is perturbed ($\delta$ is the new variance value).}
\label{fig:PLI_Unif}
\end{figure}

\begin{figure}[!ht]
  \centering
\vspace*{-3ex}	
	\includegraphics[trim=0 0 0 0,width=0.49\textwidth]{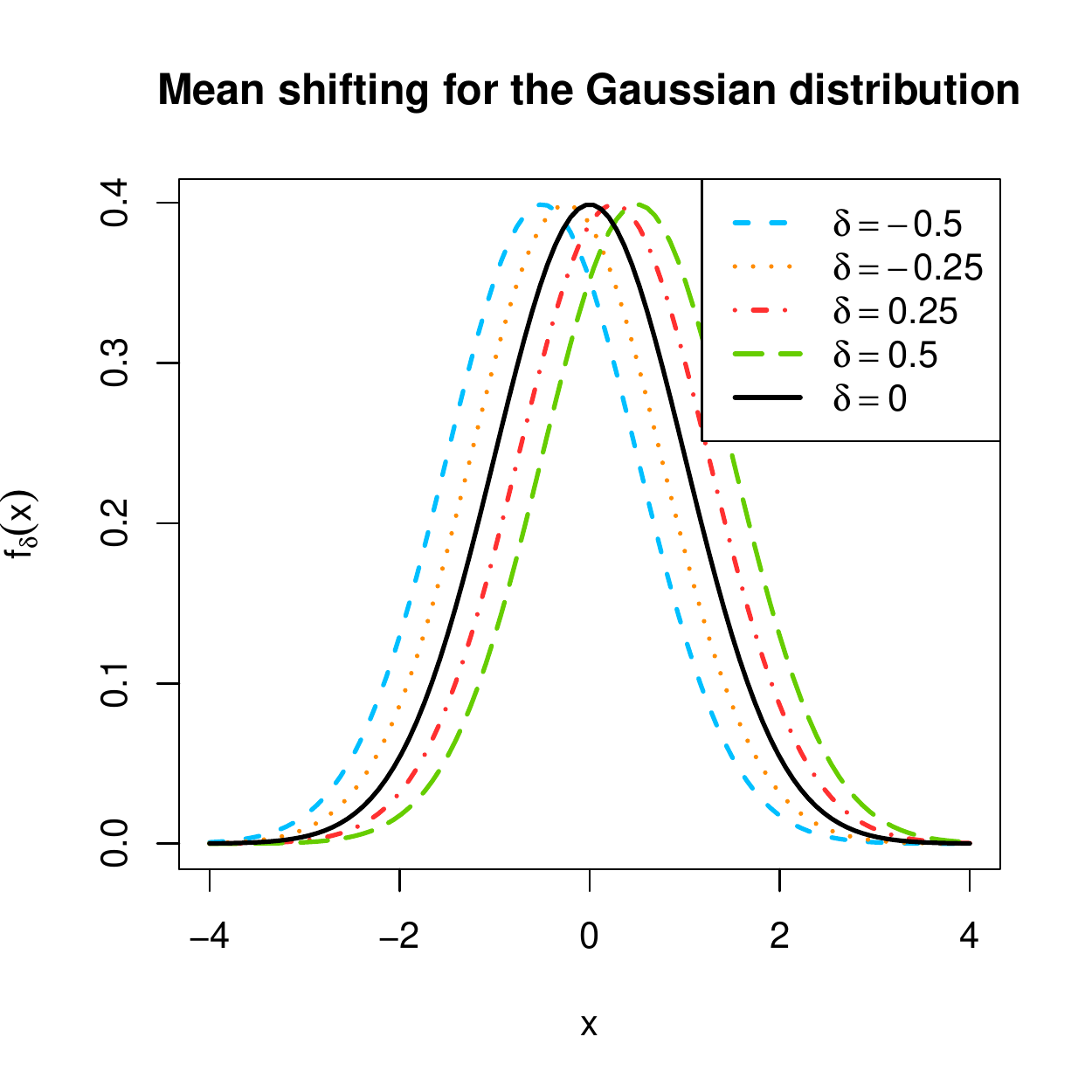}
		\includegraphics[trim=0 0 0 0,width=0.49\textwidth]{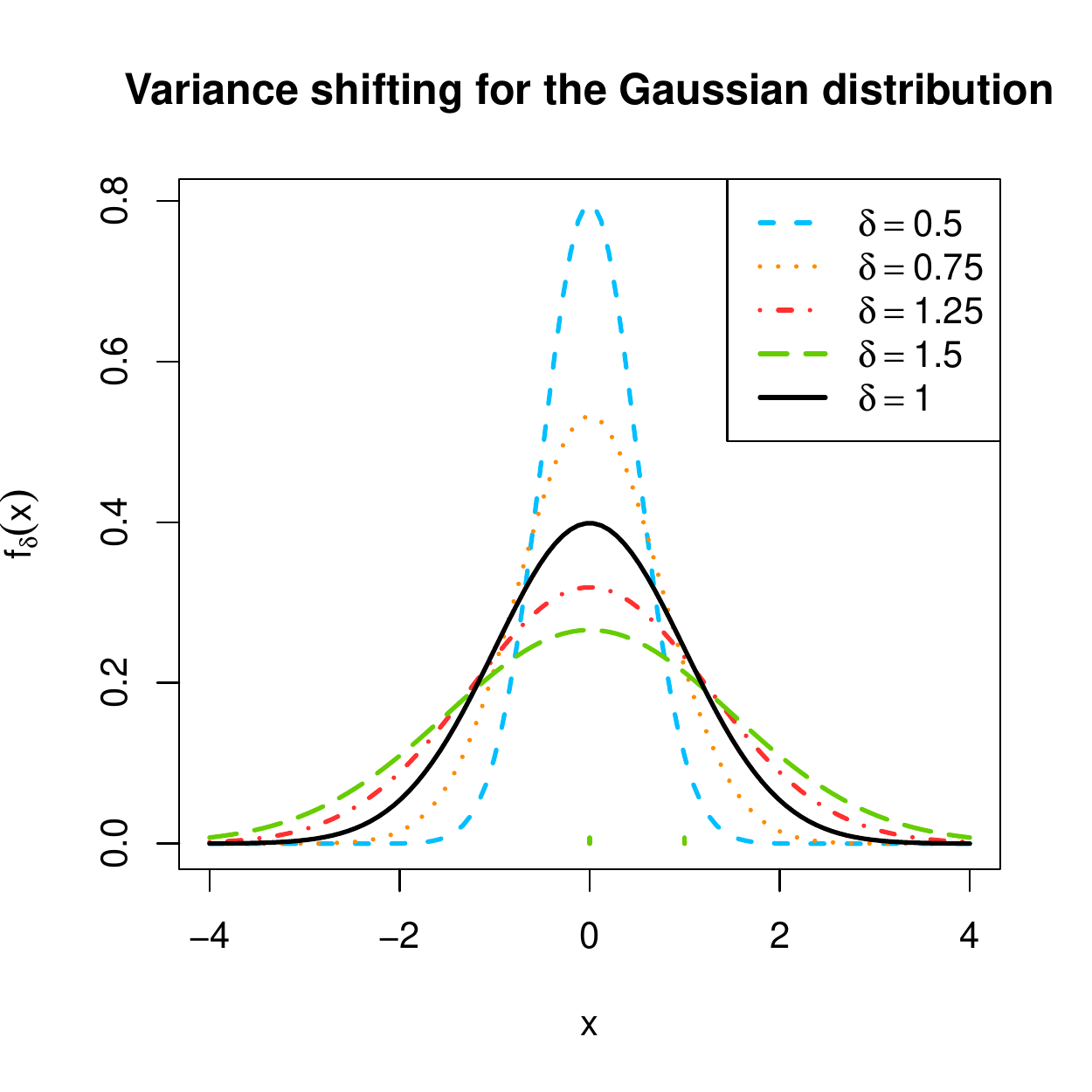}
  \caption{Examples of perturbations for the standard Gaussian pdf. 
	Left: the mean of the pdf is perturbed ($\delta$ is the new mean value). 
	Right: the standard-deviation of the pdf is perturbed ($\delta$ is the new standard-deviation value).}
\label{fig:PLI_Gauss}
\end{figure}

Figure \ref{fig:PLI_Unif} gives some examples of perturbations for the uniform pdf on $[0,1]$ (mean is $0.5$ and variance is $1/12$).
Figure \ref{fig:PLI_Gauss} gives some examples of perturbations for the standard Gaussian pdf (mean is $0$ and variance is $1$).

\begin{figure}[!ht]
  \centering
\vspace*{-3ex}	
	\includegraphics[trim=0 0 0 0,width=0.7\textwidth]{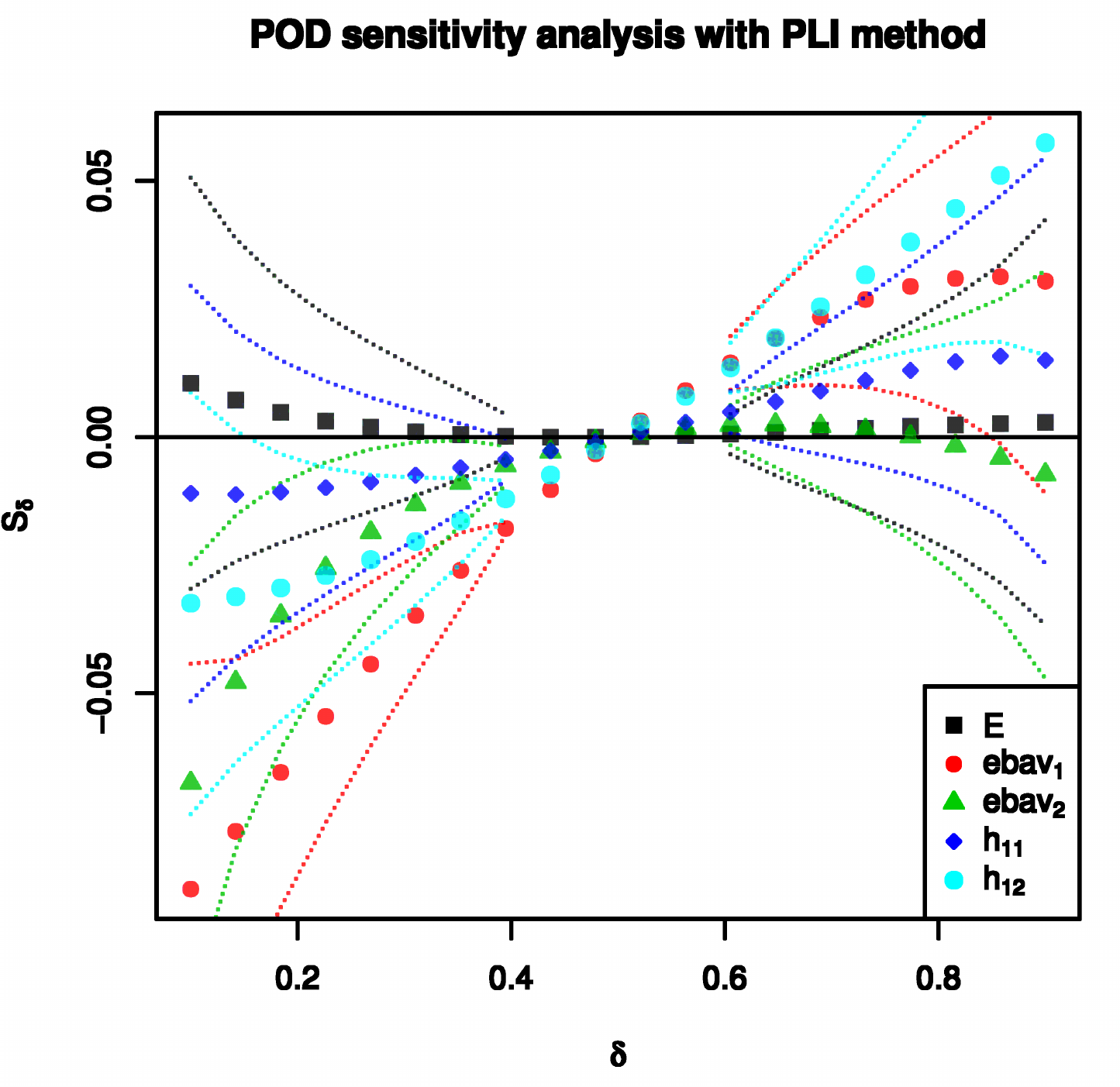}
  \caption{Graphical example for the PLI measures on the POD curves. Each color corresponds to one input parameter. The dots show the PLI values in function of the perturbation $\delta$ on the input pdf mean, while the dotted curves correspond to the $95\%$ confidence intervals on the PLI estimates.}
  \label{fig:PLI_Gauss_exemple}
\end{figure}

Figure \ref{fig:PLI_Gauss_exemple} gives an example of a PLI-based sensitivity analysis on our NDT test case (see section \ref{sec:POD}) which aims to estimate POD curves and the associated sensitivity indices to its five physical input parameters (called E, ebav$_1$, ebav$_2$, h$_{11}$ and h$_{12}$), the defect size $a$ being fixed here at a given value.
The input parameters pdfs are all uniform on $[0,1]$. 
Their non-perturbed mean is then $1/2$. 
Then, the mean of each input is modified in the range of $\delta \in [0.1,0.9]$. 
For each $\delta$ value and each influential parameter, Figure \ref{fig:PLI_Gauss_exemple} shows the PLI estimates $\hat S_{\delta}$. 
Large absolute values of $\hat S_{\delta}$ imply a large impact of the perturbation on the FRC.
Moreover, the sign of $\hat S_{\delta}$ indicates if the new probability has decreased (negative case) or increased (positive case). 
Confidence intervals shown in Figure \ref{fig:PLI_Gauss_exemple} are obtained using a Gaussian asymptotic property of the estimates $\hat S_{i,\delta}$ \cite{lemser13}. 
In Figure  \ref{fig:PLI_Gauss_exemple}, one can see for instance that increasing the h$_{12}$ input mean largely increases the FRC, while decreasing the ebav$_1$ input mean largely decreases the FRC.
In contrary, perturbations of the $E$ mean have no effect on the FRC.

\begin{figure*}[!ht]
  \centering
\vspace*{-3ex}	
	\includegraphics[width=0.49\textwidth]{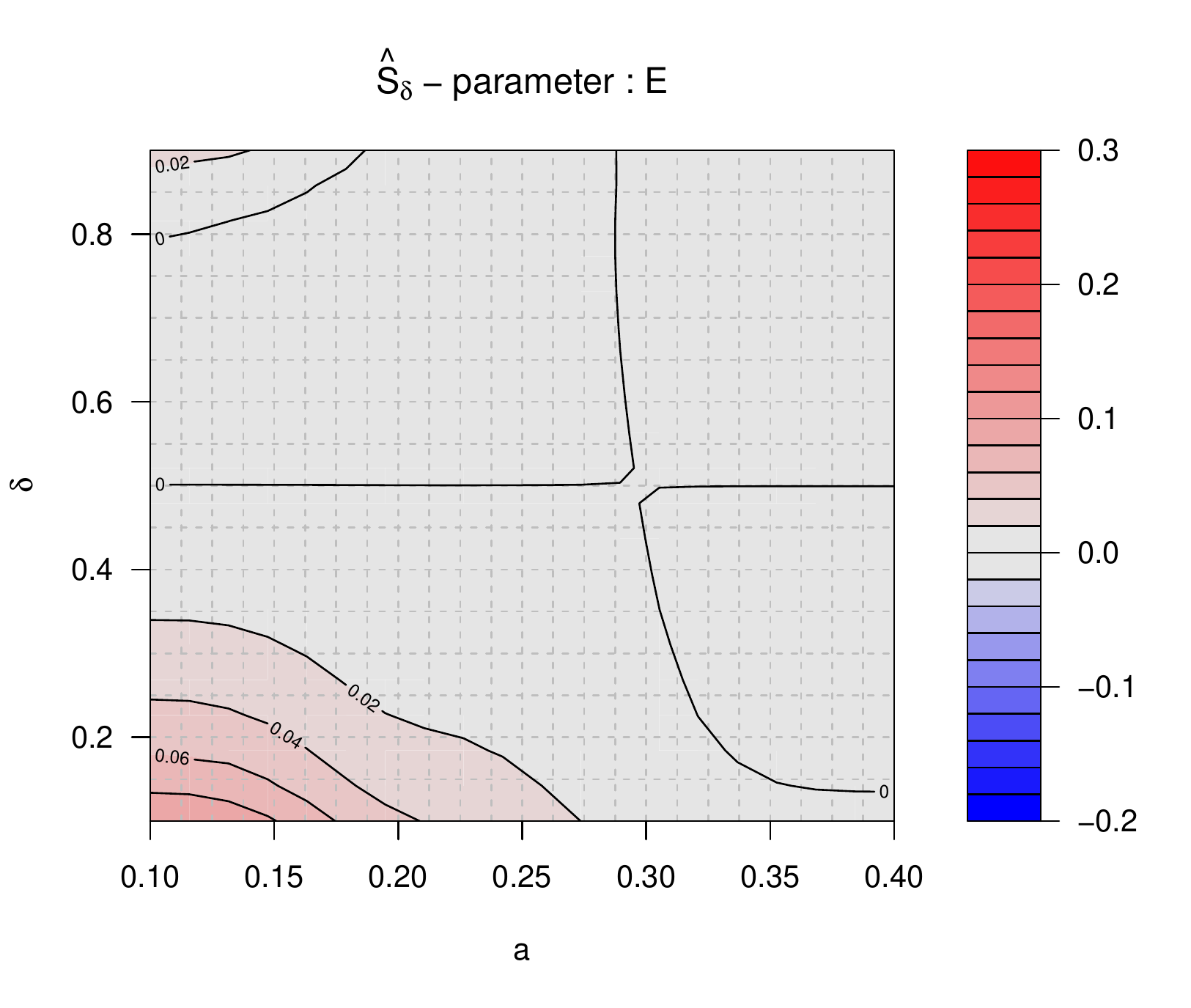}	\includegraphics[width=0.49\textwidth]{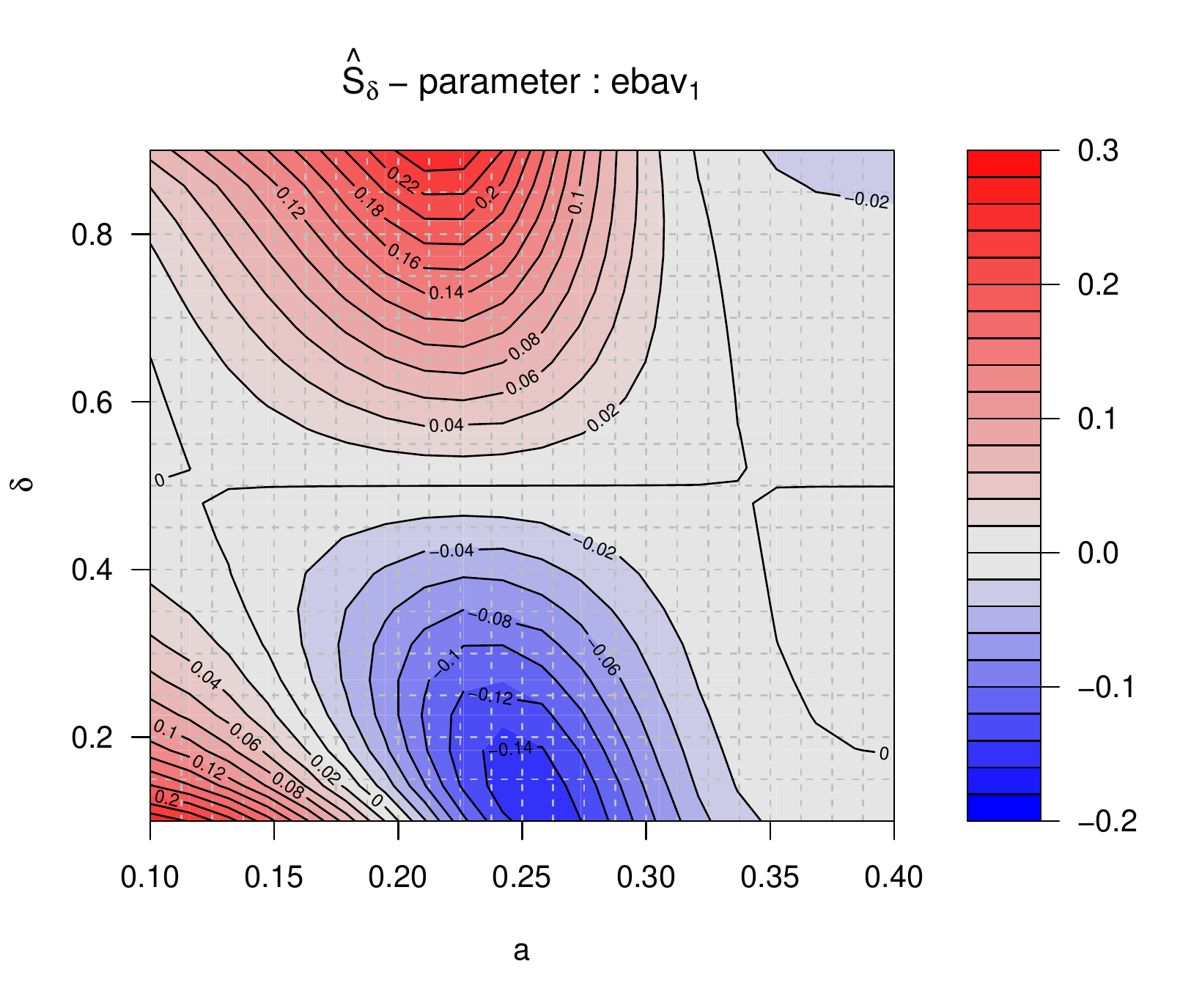}
		\includegraphics[width=0.49\textwidth]{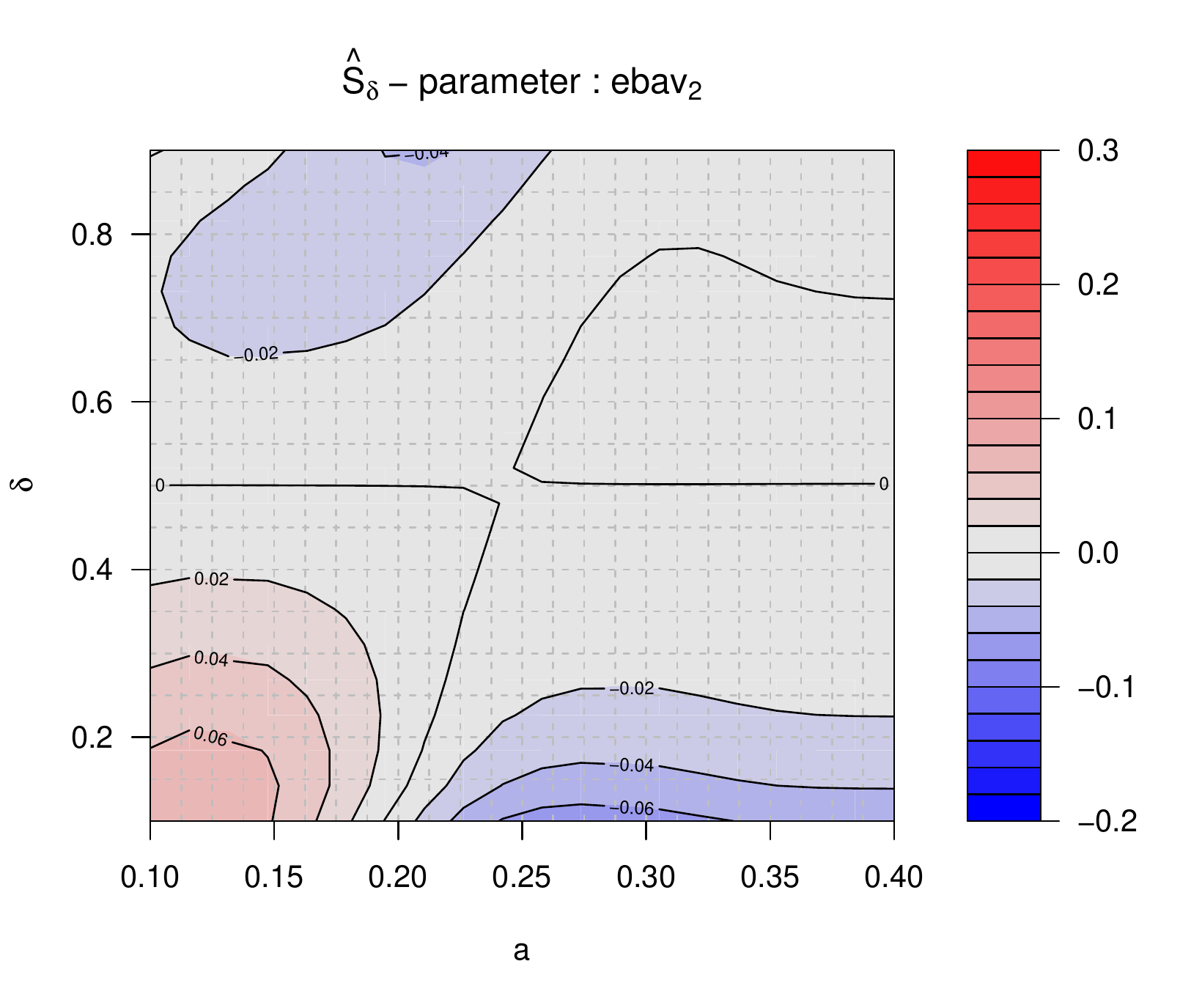}	\includegraphics[width=0.49\textwidth]{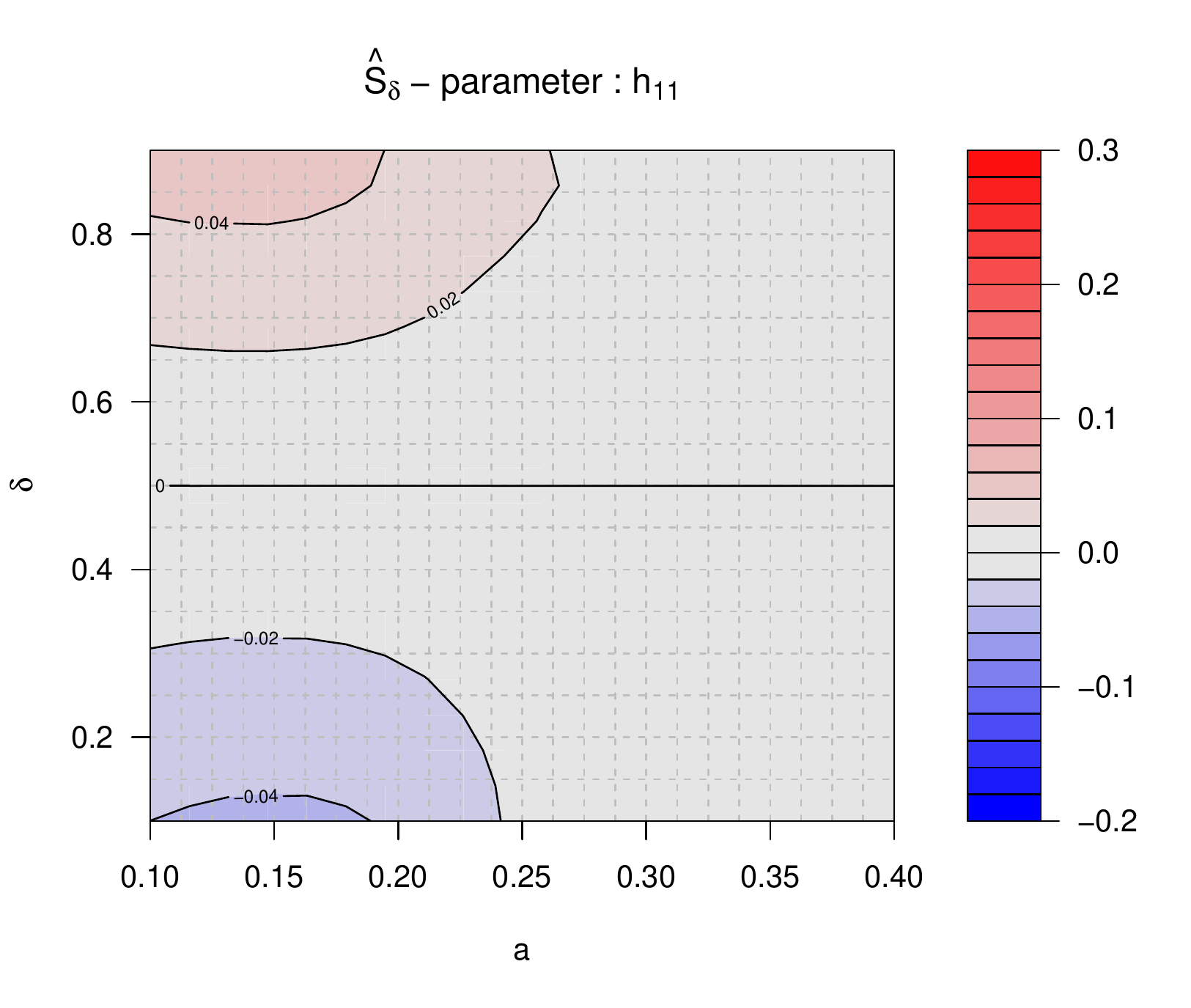}
\includegraphics[width=0.49\textwidth]{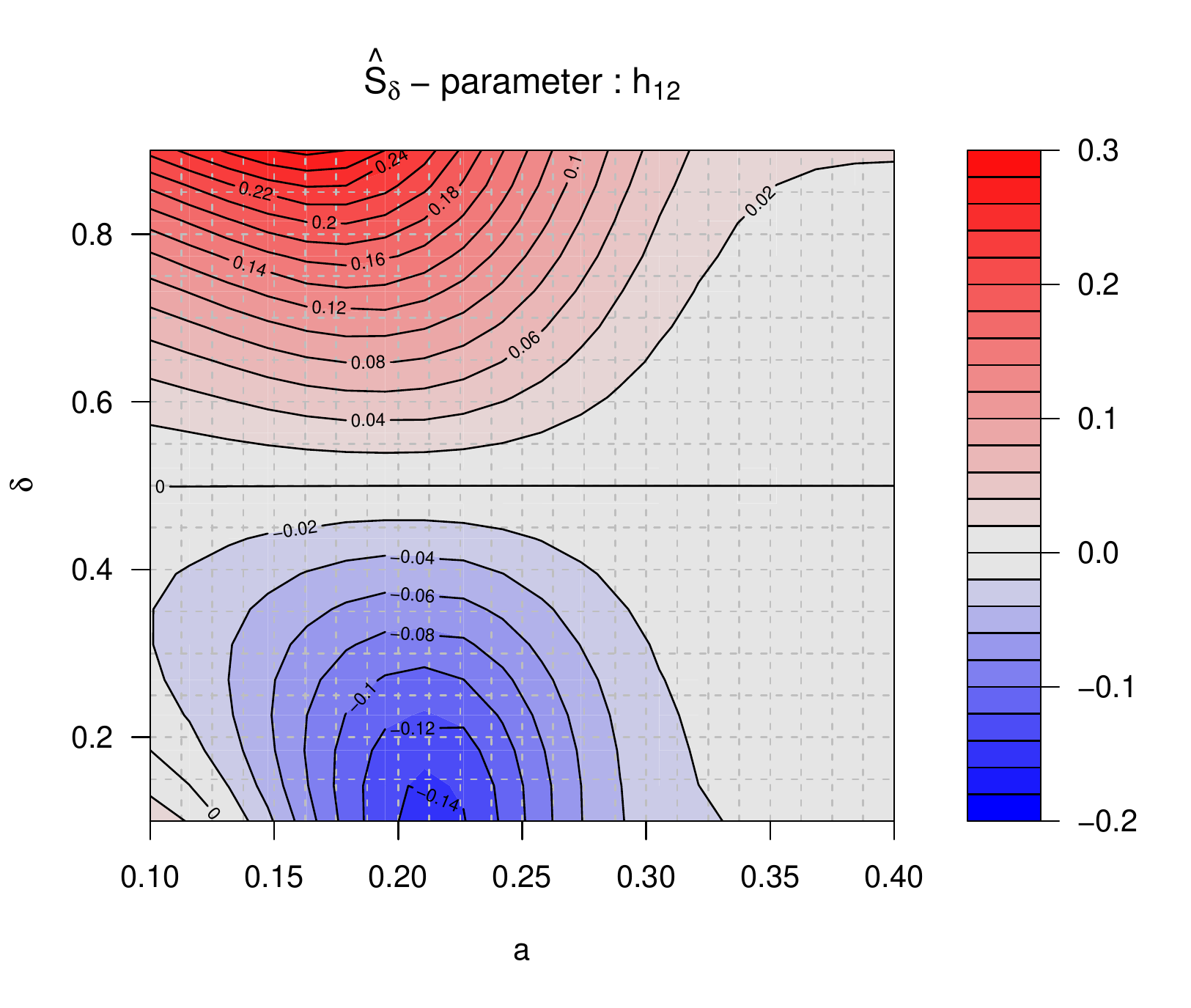}
  \caption{PLI measures on the POD curves for different perturbations on the pdf means and different values of $a$. 
	$\delta=0.5$ corresponds to the non perturbed case. 
	The grey color corresponds to no change in the FRC, while the red color (resp. blue) corresponds to an increase (resp. decrease) of the FRC.}
  \label{fig:PLI03a}
\end{figure*}

In order to have a global view of the FRC sensitivity for different defect sizes $a$, Figure \ref{fig:PLI03a} gives the PLI graphs for each input, with variations of $a$ and the $\delta$ perturbation.
Clearly, the inputs ebav$_1$ and h$_{12}$ have strong impacts on the FRC, but especially in the $a<0.3$ range.
In contrary, the inputs $E$ and h$_{11}$ have no influence. 
These elements are key points for the engineers to understand the effects of the physical parameters on the POD curve.

\section{Conclusion}

In this paper, several methods of uncertainty and sensitivity analysis of model outputs have been developed on the particular quantity of interest of the functional risk curve (FRC).
A FRC provides the probability of an undesirable event as a function of the value of a critical input parameter of a considered physical system.
Focus has been on the use of the Gaussian process metamodels in order to build FRCs from numerical experiments.
This approach is useful when the computer model is expensive to evaluate in such a way that only a small sample of the model output can be obtained.
In addition to the mean risk curve, the metamodel allows to obtain the confidence bands via conditional Gaussian process simulations.

One important advantage of using a metamodel is that sensitivity analysis is facilitated.
We have defined two kinds of sensitivity indices related to the FRC as the quantity of interest.
First, we have formulated the FRC aggregated Sobol' indices, which are variance-based measures.
Second, based on perturbation of the pdf of each model input variable, sensitivity indices (called PLI) of the model inputs on the FRC are also proposed.
These allow to understand the effect of misjudgment on the pdf of each input parameters.

An example taken from simulated NDE inspections highlights the added value of the FRC, which in this context correspond to POD curves. The two proposed sensitivity analysis methods provide strongly complementary information to experts interpreting the NDE results.
However, FRCs are used in many other engineering contexts, e.g. in the seismic fragility assessment (as shown in section \ref{sec:seismic}) and in the evaluation of hydraulic works reliability subject to extremely high water levels \cite{pas14}.
These tools would also be useful in separating the effects of stochastic input variables and epistemic input parameters.
For instance, this problem is highlighted in \cite{heljoh10} and \cite{paskel12}.

Finally, two main mathematical perspectives are identified from this work.
First, adaptive designs can be developed: in the spirit of \cite{becgin12}, it would consist of using the Gaussian process model in order to define some SUR (``Sequential Uncertainty Reduction'') criteria on the FRC as the quantity of interest.
Optimizing such a criterion would provide a new set of input parameter values, which would be run with the computer code in order to decrease the POD confidence interval. 
Second, \cite{brofor15} have started a preliminary study that considers FRC as a random distribution function.
This framework allows to deal with stochastic computer codes instead of the deterministic ones of this paper (see for example Kleijnen \cite{kle15}).

\section*{Acknowledgments}

Part of this work has been backed by French National Research Agency (ANR) through project ByPASS ANR-13-MONU-0011.
We acknowledge three anonymous reviewers for their deep analysis of the paper and the useful comments they provide.
We are also grateful to Alberto Pasanisi for giving helpful discussions and Andreas Schumm for his help with the English language.

\section*{References}
\bibliographystyle{elsarticle-num}

\end{document}